\documentclass{amsart}
\usepackage[all]{xy}
\usepackage{amscd,amssymb}
\usepackage{graphicx}
\usepackage{hyperref}
\hypersetup{
    colorlinks=true,       
    linkcolor=blue,          
    citecolor=blue,        
    filecolor=blue,      
    urlcolor=blue           
}
\usepackage[alphabetic,backrefs,lite]{amsrefs}

\newcommand{\defi}[1]{\textsf{#1}}

\newtheorem{theorem}{Theorem}[section]
\newtheorem{lemma}[theorem]{Lemma}
\newtheorem{prop}[theorem]{Proposition}

\newtheorem{cor}[theorem]{Corollary}

\theoremstyle{definition}
\newtheorem{definition}[theorem]{Definition}
\newtheorem{example}[theorem]{Example}

\theoremstyle{remark}
\newtheorem{remark}[theorem]{Remark}

         \newcommand{\cO}{\mathcal O}

\newcommand{\bK}{\mathbb K}

\newcommand{\bP}{\mathbb P}
\newcommand{\bQ}{\mathbb Q}
\newcommand{\bZ}{\mathbb Z}

\newcommand{\Nef}{{\rm Nef}}
\newcommand{\Aut}{{\rm Aut}}
\newcommand{\Pic}{{\rm Pic}}

\newcommand{\Eff}{{\rm Eff}}

\newcommand{\lE}{\langle E \rangle}
\newcommand{\SAmple}{{\rm SAmple}}

\begin{document}

\title{Nef and semiample divisors on rational surfaces}
\thanks{The first author was supported by Proyecto FONDECYT Regular 2011, N. 1110096. The second author was supported by EPSRC grant number EP/F060661/1.}
\author{Antonio Laface}
\address{
Departamento de Matem\'atica, 
Universidad de Concepci\'on, 
Casilla 160-C,
Concepci\'on, Chile}
\email{antonio.laface@gmail.com}

\author{Damiano Testa}

\address{Mathematical Institute, 24-29 St Giles', Oxford OX1 3LB, United Kingdom}
\email{adomani@gmail.com}

\date{\today}

\begin{abstract}
In this paper we study smooth projective rational surfaces, defined over an algebraically closed field of any characteristic, with pseudo-effective anticanonical divisor.  We provide a necessary and sufficient condition in order for any nef divisor to be semiample.  We adopt our criterion to investigate Mori dream surfaces in the complex case.
\end{abstract}

\maketitle

\section*{Introduction}

Let $X$ be a smooth projective rational surface defined over an algebraically closed field $\bK$ of any characteristic.  A problem that has recently attracted attention consists in finding equivalent characterizations of Mori dream surfaces, that is surfaces whose Cox ring is finitely generated (see~\cite{adhl} for basic definitions), for instance in terms of the Iitaka dimension of the anticanonical divisor. Indeed, if the Iitaka dimension of the anticanonical divisor is $2$, then $X$ is always a Mori dream surface, as shown in~\cite{tvv}. If the Iitaka dimension of the anticanonical divisor is $1$, then $X$ admits an elliptic fibration $\pi \colon X\to \bP^1$ and it is a Mori dream surface if and only if the relatively minimal elliptic fibration of $\pi$ has a jacobian fibration with finite Mordell-Weil group, as shown in~\cite{al}. The authors are not aware of any previously known example of a Mori dream surface with Iitaka dimension of the anticanonical divisor equal to $0$.

By~\cite{hk}*{Definition~1.10}, a smooth projective surface $X$ is a Mori dream space if the irregularity of $X$ vanishes and if the nef cone of $X$ is the affine hull of finitely many  semiample classes.  In this paper we concentrate on giving necessary and sufficient conditions in order for any nef divisor to be semiample in the case in which $X$ is a smooth projective rational surface with non-negative anticanonical Iitaka dimension.  Our theorems hold in any characteristic; in the complex case we provide examples of Mori dream surfaces with vanishing anticanonical Iitaka dimension.

The paper is organized as follows. We introduce necessary notation in Section~\ref{nota}.  In Section~\ref{nesa} we prove some technical lemmas about the stable base locus of divisors on a smooth projective rational surface. We use these lemmas to prove Corollary~\ref{ii} of Section~\ref{kneg}.  In Section~\ref{knoteff} we prove that, if the anticanonical divisor is not effective, then every nef divisor is semiample. In Section~\ref{sekeff} we consider the case where the anticanonical divisor is effective, providing a necessary and sufficient condition for every nef divisor to be semiample.  At the end of this section we prove that the blow-up at 10 or more distinct points on a smooth cubic curve is never a Mori dream space. Finally, in Section~\ref{seexa} we give an example of a Mori dream surface with vanishing anticanonical Iitaka dimension.

\subsection*{Ackowledgments}
It is a pleasure to thank Igor Dolgachev for the suggestion of looking at Coble surfaces with finite automorphism group to produce examples of Mori dream surfaces, and to thank Michela Artebani for helpful remarks and discussions. We would also like to thank Ana-Maria Castravet for pointing out a mistake in an older version of this paper.
We have made large use of the computer algebra program Magma~\cite{Magma} in Section~\ref{seexa} of this paper.

\section{Notation} \label{nota}

In this section we gather the basic definitions and the standard conventions that we use throughout the paper.  
Let $X$ be a smooth projective surface with ${\rm H}^1(X,\cO_X)=0$. 
All the cones mentioned in this paper are contained in the 
rational vector space $\Pic_\bQ(X) = \Pic(X)\otimes_{\bZ}\bQ$, see~\cite{adhl}.
The {\emph{effective cone}} $\Eff(X)$ consists of the numerical equivalence classes of divisors 
admitting a positive multiple numerically equivalent to an effective divisor; the {\emph{pseudo-effective cone}} $\overline{\Eff}(X)$ is the closure of the effective cone.
The {\emph{nef cone}} $\Nef(X)$ consists of the classes of divisors $D$ such that
$D\cdot C\geq 0$ for any curve $C$ in $X$. The pseudo-effective
cone and the nef cone are dual with respect to each other under the pairing induced by the intersection
form of $\Pic(X)$.
The {\emph{semiample cone}} $\SAmple(X)$ consists of the classes of 
divisors admitting a positive multiple whose associated linear system
is base point free.
Clearly, every semiample divisor is also nef, and therefore there is an inclusion $\SAmple(X) \subset \Nef(X)$.

We specialize Definition~1.10 of~\cite{hk} to our context.

\begin{definition}
Let $X$ be a smooth projective surface. We say that $X$ is
a {\em Mori dream surface} if the following conditions hold:
\begin{enumerate}
\item 
the group ${\rm H}^1(X,\cO_X)$ vanishes;
\item \label{fine}
the nef cone of $X$ is generated by a finite number
of semiample classes.
\end{enumerate}
\end{definition}

Condition~\eqref{fine} means that the cones $\Nef(X)$ and $\SAmple(X)$ coincide and that both
are polyhedral. The statement that the cone $\Eff(X)$ is polyhedral and that the cones $\Nef(X)$ and $\SAmple(X)$ coincide is also equivalent to Condition~\eqref{fine}.
If the cone $\Eff(X)$ is polyhedral and ${\rm rk}(\Pic(X)) \geq 3$, it is not difficult to prove 
that $\Eff(X)$ is spanned by classes of 
{\emph{negative curves}}, that is integral curves $C$ 
with $C^2<0$ (see~\cite{al}*{Proposition~1.1}). Examples of negative curves
are {\emph{$(-n)$-curves}} which are smooth
rational curves of self-intersection $-n$.
An {\emph{exceptional curve}} is an alternative name for a $(-1)$-curve.

Recall that the {\emph{Iitaka dimension}} $k(D)$ of a divisor $D$ is
one less than the Krull dimension of the ring 
\[
\bigoplus _{n \geq 0} {\rm H}^0(X,nD)
\]
if a positive multiple of $D$ is effective and $-\infty$ otherwise.
Our initial approach is to study Mori dream surfaces according
to their anticanonical Iitaka dimension $k(-K_X)$.
Toric surfaces and {\emph{generalized del Pezzo surfaces}},
that is rational surfaces with nef and big anticanonical divisor, have
anticanonical Iitaka dimension $2$.
Rational elliptic surfaces have anticanonical Iitaka dimension $1$.

\section{Nef and semiample divisors} \label{nesa}

Let $X$ be a smooth projective surface and let $C$ be an effective divisor on $X$.  Suppose that the matrix of pairwise intersections of the components of $C$ is negative definite.  It follows that there is an effective divisor $B$ with support on $C$ such that for every irreducible component $D$ of $C$ we have $B \cdot D < 0$.  

\begin{definition}
We say that an effective divisor $B$ with support in $C$ is a \defi{block on $C$} if for every irreducible component $D$ of $C$ the inequality $B \cdot D < \min \{0,-K_X \cdot D\}$ holds.
\end{definition}

Often in our applications, the divisor $C$ is the union of all the integral curves of $X$ orthogonal to a big and nef divisor $N$ on $X$; in this case the required negative definiteness follows by the Hodge Index Theorem.  In such cases we say that a block on $C$ is a \defi{block for $N$}.
Observe that from the definition of block it follows easily that if $B$ is a block for $C$, then the supports of $B$ and of $C$ coincide.

\begin{lemma} \label{naka}
Let $X$ be a smooth projective surface.  Let $N$ be a big and nef divisor on $X$ and let $B$ 
be a block for $N$.  Then the stable base locus of $N$ is the stable base locus of the line bundle $\cO_B(N)$.  
In particular, if $R$ is a connected component of the block $B$ such that the group ${\rm H}^1(R,\cO_R)$ 
vanishes, then $R$ is disjoint from the stable base locus.
\end{lemma}

\begin{proof}
To prove the lemma, we use some of the techniques of the proof of~\cite{Nak}*{Theorem~1.1}.  
Let $E$ be an effective divisor on $X$; we establish below the following assertion: 
\renewcommand{\labelitemi}{($\star$)}
\begin{itemize}
\item
there is an effective divisor $Z$ with the same support as $B \cup E$ such that, for sufficiently large $n$, the restriction map 
\[
{\rm H}^0 \bigl( X, \mathcal{O}_X (nN) \bigr) \longrightarrow {\rm H}^0 \bigl( X, \mathcal{O}_{Z} (nN) \bigr)
\]
is surjective.
\end{itemize}
Assuming ($\star$) we deduce that, if $E$ is an integral curve with $N \cdot E > 0$, then the restriction of $nN$ to $Z$ admits global sections not vanishing identically on $E$, and hence that $E$ is not contained in the stable base locus of $N$.  Therefore the support of the stable base locus of $N$ is contained in the support of $B$.  Finally, from the proof of ($\star$) we obtain that when $E=0$ we can choose $Z=B$, establishing the main statement of the lemma.

We now turn to the proof of ($\star$).  Replacing if necessary 
$B$ by a positive multiple, we obtain a block $B'$ for $N$ such that for every irreducible component $D$ of $B$ we 
have $B' \cdot D < -(E + K_X) \cdot D$; note that, in the case in which $E=0$, we may choose $B'=B$.  We show that the divisor $Z = B'+E$ satisfies ($\star$).  Let $A$ be an ample divisor on $X$.  Since the divisor $N$ 
is big, for large enough $n$ the divisor $F=nN-Z-K_X-A$ is effective and hence the divisor $nN-Z-K_X = A+F$ 
is big.  Moreover, if $D$ is an integral curve on $X$ such that $D \cdot (nN-Z-K_X) \leq 0$, then $D$ is a component 
of $F$ and it is not a component of $B'$.  It follows that the inequality $N \cdot D>0$ holds, and since $F$ has only 
finitely many components, we obtain that for large enough $n$ the divisor $nN-Z-K_X$ has positive intersection 
with every curve on $X$; since it is also big, we deduce from the Nakai-Moishezon criterion that the divisor $nN-Z-K_X$ is 
ample.

Thus, for every sufficiently large integer $n$ the Kodaira Vanishing Theorem 
implies that the group ${\rm H}^1 (X, \cO_X(nN-Z))$ vanishes and ($\star$) follows from the long exact cohomology sequence associated to the sequence 
\[
0 \longrightarrow \mathcal{O}_X(nN-Z) \longrightarrow \mathcal{O}_X (nN) \longrightarrow 
\mathcal{O}_{Z} (nN) \longrightarrow 0 ,
\]
establishing ($\star$).

For the last assertion, by the assumptions it follows that the restriction of $N$ to $R$ is trivial; in particular the line bundle $\cO_R(N)$ is globally generated and the result follows from what we already proved.
\end{proof}

\begin{lemma} \label{ga}
Let $X$ be a smooth projective surface with $q(X) = 0$ and let $B$ be an effective divisor on $X$.  If the divisor 
$K_X+B$ is not effective, then the group ${\rm H}^1(B,\cO_B)$ vanishes.  In particular, if there is a nef divisor 
$N$ such that the inequality $N \cdot (K_X + B) < 0$ holds, then the group ${\rm H}^1(B,\cO_B)$ vanishes.
\end{lemma}

\begin{proof}
From the sequence 
\[
0 \longrightarrow \mathcal{O}_X(-B) \longrightarrow \mathcal{O}_X \longrightarrow \mathcal{O}_B \longrightarrow 0
\]
we deduce that it suffices to show that the group ${\rm H}^2(X,\cO_X(-B))$ vanishes.  The first part follows by Serre duality.  
The second part is an immediate special case of the first and the lemma follows.
\end{proof}

\section{The $K_X$-negative part of the nef cone} \label{kneg}

\begin{lemma} \label{bbp}
Let $X$ be a smooth projective rational surface.  If the divisor $N$ on $X$ is nef and it satisfies $-K_X \cdot N > 0$, 
then $N$ is semiample.  If $N$ is not big, then it is base point free.
\end{lemma}

\begin{proof}
Suppose first that $N$ is big and let $B$ be a block for $N$.  The result follows by Lemmas~\ref{naka} and~\ref{ga} since the 
inequality $N \cdot (K_X+B) < 0$ holds.

Suppose now that $N$ is not big, and hence that $N^2=0$.  The divisors $K_X-N$ and $K_X+N$ are not effective since 
$N \cdot (K_X-N)<0$ and $N \cdot (K_X+N)<0$.  Therefore, by the Riemann-Roch formula, we deduce that the divisor $N$ 
is linearly equivalent to an effective divisor.  To avoid introducing more notation, we assume that $N$ itself is effective.  We 
deduce that the base locus of $N$ is contained in $N$ and that the group ${\rm H}^1(N,\cO_N)$ vanishes by 
Lemma~\ref{ga}.  In particular the line bundle $\cO_N(N)$ is trivial and from the exact sequence 
\[
0 \longrightarrow \mathcal{O}_X \longrightarrow \mathcal{O}_X (N) \longrightarrow \mathcal{O}_N (N) 
\longrightarrow 0
\]
and the vanishing of the group ${\rm H}^1(X,\cO_X)$, we deduce that the base locus of $|N|$ is empty.
\end{proof}

As a consequence of what we proved so far, we deduce the following result (see also~\cite{al}*{Theorem~3.4} and~\cite{tvv}).

\begin{cor} \label{ii}
Let $X$ be a smooth projective rational surface and suppose that there are a positive integer $k$ and $\bQ$-divisors $P$ 
and $E$ such that $P$ is non-zero and nef, $E$ is effective, $-K_X=P+E$ and $kP$ is integral.  Every nef divisor on $X$
not proportional to $P$ is semiample; 
moreover $P$ itself is semiample if $P$ is big or $kP$ is base point 
free.
\end{cor}

\begin{proof}
Let $N$ be a non-zero nef divisor on $X$.  If $-K_X \cdot N > 0$, then we conclude that $N$ is semiample by 
Lemma~\ref{bbp}.  Suppose that $-K_X \cdot N = 0$.  In particular we deduce that $P \cdot N = 0$ and hence 
that $P$ and $N$ are not big and that $N$ is proportional to $P$.  Again Lemma~\ref{bbp} allows us to conclude.
\end{proof}

\section{The anticanonical divisor is not effective} \label{knoteff}

Let $X$ be a smooth projective rational surface.  In this section we assume that the anticanonical divisor $-K_X$ is not linearly 
equivalent to an effective divisor, but a positive multiple of $-K_X$ is.  Let $e$ be the least positive integer such 
that the linear system $|{-eK_X}|$ is not empty and let $E$ be an element of $|{-eK_X}|$.

\begin{theorem} \label{nof}
Every nef divisor on $X$ is semiample.
\end{theorem}

\begin{proof}
If $N$ is trivial, then the result is clear; we therefore suppose that 
$N$ is non-trivial.  Note that the divisor $K_X-N$ is not effective, as it is the negative of a non-zero pseudo-effective 
divisor.  Thus by the Riemann-Roch formula and Serre duality we deduce that $N$ is linearly equivalent to 
an effective divisor; to avoid introducing more notation, we assume that $N$ itself is effective.  If $N$ is big, then 
let $B$ denote a block for $N$; if $N$ is not big, then let $B=N$.  In both cases the base locus of $N$ is contained 
in $B$.

Suppose first that $h^1(B,\cO_B) = 0$.  By Lemma~\ref{naka} we deduce that if $N$ is big, then it is semiample.  
Thus we reduce to the case in which $N$ is not big; the line bundle $\cO_N(N)$, having non-negative degree 
on each irreducible component of its support, is globally generated.  From the sequence 
\[
0 \longrightarrow \mathcal{O}_X \longrightarrow \mathcal{O}_X (N) \longrightarrow \mathcal{O}_N (N) 
\longrightarrow 0
\]
and the vanishing of the group ${\rm H}^1(X,\cO_X)$ we deduce that also in this case the divisor $N$ is semiample.

Suppose now that $h^1(B,\cO_B) \geq 1$.  By Lemma~\ref{ga} we obtain that the divisor $K_X+B$ is linearly 
equivalent to the effective divisor $B'$.  We deduce that the divisor $eB$ is linearly equivalent to the divisor 
$eB'+E$.  If the effective divisor $eB$ were equal to the effective divisor $eB'+E$, it would follow that every prime 
divisor appears with multiplicity divisible by $e$ in $E$; in particular, the divisor $E$ would be a multiple of an 
effective divisor linearly equivalent to $-K_X$.  Since we are assuming that $-K_X$ is not effective, it follows that 
$eB$ and $eB'+E$ are not equal and we deduce that the dimension of $|eB|$ is at least one.  In particular, $N$ is 
not big, since otherwise the matrix of pairwise intersections of the components of $B$ would be negative definite and 
hence no multiple of $B$ would move.  Therefore we have $B=N$; write $eN=M+F$ where $M \neq 0$ is the moving 
part of $|eN|$ and $F$ is the fixed part.  Since the divisor $N$ is nef and not big it follows that $M \cdot F=F^2=0$ and 
by the Hodge Index Theorem we obtain that $F$ is proportional to $M$.  Since $M$ is semiample, the result follows.
\end{proof}

\section{The anticanonical divisor is effective} \label{sekeff}

Let $X$ be a smooth projective rational surface.  In this section we assume that $X$ admits an effective anticanonical 
divisor.  We let $E$ be an element of the linear system $|{-K_X}|$ and we denote by $\omega_E$ the dualizing sheaf 
of $E$.  Since $E$ is an effective divisor on a smooth surface, its dualizing sheaf is invertible and by the adjunction 
formula we have $\omega_E \simeq \cO_E$.  Moreover from the sequence 
\[
0 \longrightarrow \mathcal{O}_X(-E) \longrightarrow \mathcal{O}_X 
\longrightarrow \mathcal{O}_E \longrightarrow 0
\]
we deduce that $h^0(E,\cO_E)=1$; in particular the divisor $E$ is connected and has arithmetic genus 
equal to one.

Let $\iota \colon E \to X$ denote the inclusion and let 
\[
\iota ^* \colon \Pic (X) \longrightarrow \Pic (E)
\]
be the pull-back map induced by $\iota$.  The map $\iota ^*$ is a homomorphism of abelian groups whose image 
is a finitely generated subgroup of $\Pic(E)$, since $\Pic(X) \simeq \bZ^{10-(K_X)^2}$.  
If $N$ is a divisor on $X$, then we denote by $N_E \in \Pic (E)$ the class of the line bundle $\cO_E(N)$.

We let $\Gamma \subset \Pic (E)$ be the image of the lattice $\Nef(X) \cap \langle K_X \rangle^\perp$ under the homomorphism $\iota^*$; equivalently $\Gamma$ is the abelian subgroup of $\Pic(E)$ generated by the classes $N_E$, where $N$ ranges among the nef divisors such that $K_X \cdot N = 0$.  In particular, $\Gamma$ is a finitely generated abelian group.  We shall prove in Theorem~\ref{fatto} that the equality of semiample and nef divisors is equivalent to the finiteness of the group $\Gamma$.

\begin{lemma} \label{tosto}
Let $N$ be a nef divisor satisfying $N \cdot K_X = 0$.  If $N$ is semiample, then $N_E \in \Pic(E)$ is a torsion 
element.
\end{lemma}

\begin{proof}
Let $r$ be a positive integer such that the linear system $|rN|$ is base point free.  Since $N \cdot K_X = 0$ we deduce 
that the linear system $|rN|$ contains a divisor disjoint from $E$ and hence $(rN)_E$ represents the trivial line bundle 
and $N_E \in \Pic(E)[r]$.
\end{proof}

\begin{lemma} \label{bb}
Let $N$ be a big and nef divisor satisfying $N \cdot K_X = 0$.  The support of the stable base locus of $N$ 
is contained in $E$.  If moreover $N_E$ is a torsion element of $\Pic(E)$, then $N$ is semiample.
\end{lemma}

\begin{proof}
First, we reduce to the case in which all the exceptional curves $F$ on $X$ such that $N \cdot F = 0$ are 
contained in $E$.

Let $F \subset X$ be an exceptional curve such that 
$F \cdot N = 0$ and $F$ is not contained in $E$, and let $b \colon X \to X'$ be the contraction of $F$.  By 
construction, the divisor $N$ is the pull-back of a big and nef divisor $N'$ on $X'$; moreover since 
$E \cdot F = 1$ and $F$ is not contained in $E$, it follows that $F$ meets $E$ transversely at a single smooth 
point.  Therefore the restriction of $b$ to $E$ is an isomorphism to its image $E'$, and the restriction of 
$\cO_{X'}(N')$ to $E'$ is isomorphic to the restriction of $\cO_X(N)$ to $E$.  Finally, if $N'$ is semiample, then 
also $N$ is semiample.  Repeating if necessary the above construction starting from $X'$, $N'$ and $E'$ we 
conclude the reduction step.  Observe that the process described above terminates, since at each stage we 
contract an exceptional curve, thereby reducing the rank of the Picard group of $X$.

Thus, we assume from now on that the exceptional curves $F$ on $X$ such that $F \cdot N = 0$ are contained in $E$.

Let $B$ be a block for $N$.  Let $F$ be a component of $B$ not contained in $E$.  
Thus we have $F^2<0$ and $F \cdot E \geq 0$; it follows from the adjunction formula that $F$ is a smooth rational 
curve and either $F \cdot K_X = -1$ and $F$ is an exceptional curve, or $F \cdot K_X = 0$ and $F$ is a $(-2)$-curve 
disjoint from $E$.  By our reduction we conclude that every irreducible component of $B$ not contained in $E$ is 
disjoint from $E$; we thus have $B=B'+E'$ where the support of $E'$ is contained in the support of $E$ and the support 
of $B'$ is the union of the components of $B$ disjoint from $E$.

We now show that no component of $B'$ is contained in the stable base locus of $N$; for this it suffices to show that the 
divisor $K_X+B'$ is not effective by Lemmas~\ref{naka} and~\ref{ga}.  Since the matrix of pairwise intersections of the 
components of $B'$ is negative definite, the divisor $B'$ is the unique element of the linear system $|B'|$.  If the divisor 
$B'-E$ were linearly equivalent to an effective divisor $B''$, then the linear system $|B'|$ would contain the 
divisor $B''+E$ contradicting the fact that $|B'| = \{B'\}$.  Thus the divisor $K_X+B'$ is not effective, as we wanted to show, 
and the stable base locus of $N$ is contained in $E'$.

To prove the second part, we show that no component of $E$ is contained in the stable base locus of $N$.  
Let $r$ be a positive integer such that $(rN)_E$ is trivial.  From the 
Kawamata-Viehweg Vanishing Theorem (see~\cite{Xi}*{Corollary~1.4} for the positive characteristic case) we deduce that the group ${\rm H}^1 (X , \cO_X(rN-E))$ vanishes and therefore using the exact sequence 
\[
0 \longrightarrow \mathcal{O}_X(rN-E) \longrightarrow \mathcal{O}_X (rN) \longrightarrow \mathcal{O}_E(rN) \longrightarrow 0
\]
we deduce that there are sections of $\cO_X(rN)$ that are disjoint from $E$ and the lemma is proved.
\end{proof}

\begin{lemma} \label{nb}
Let $N$ be a nef non big divisor satisfying $N \cdot K_X = 0$.  If $N_E$ is a torsion element of $\Pic(E)$, then 
$N$ is semiample.
\end{lemma}

\begin{proof}
If the divisor $N$ is trivial, then the result is clear.  Suppose that $N$ is non-zero and let $r$ be a 
positive integer such that $(rN)_E$ is trivial.  Observe that $h^2(X,\cO_X(rN-E))=h^0(X,\cO_X(-rN)) = 0$, 
since for every ample divisor $A$ on $X$ we have $A \cdot N>0$ by Kleiman's criterion and therefore 
the divisor $-rN$ cannot be effective.  Thus from the exact sequence 
\[
0 \longrightarrow \mathcal{O}_X(rN-E) \longrightarrow \mathcal{O}_X (rN)
\longrightarrow \mathcal{O}_E(rN) \longrightarrow 0
\]
and the fact that $h^1(E,\cO_E)=1$ we deduce that $h^1(X,\cO_X (rN)) \geq 1$ and thus that 
the dimension of the linear system $|N|$ is at least one, by the Riemann-Roch formula.  We may 
therefore write $N=M+F$, where $M \neq 0$ is the moving part and $F$ is the divisorial base locus of 
$|N|$.  Since $N$ is nef and not big we deduce that $M \cdot F=F^2=0$ and thus by the Hodge 
Index Theorem we conclude that $F$ is proportional to $M$ and the lemma follows.
\end{proof}

\begin{theorem} \label{fatto}
Let $X$ be a smooth projective rational surface and suppose that $E$ is an effective divisor in the 
anticanonical linear system.  The subgroup $\Gamma$ of $\Pic(E)$ is finite if and only if every 
nef divisor on $X$ is semiample.
\end{theorem}

\begin{proof}
Suppose that the group $\Gamma $ is finite and let $N$ be a nef divisor on $X$.  If $-K_X \cdot N > 0$, 
then the result follows from Lemma~\ref{bbp}.  Suppose therefore that $-K_X \cdot N = 0$.  It follows that 
$N_E$ is contained in $\Gamma $ and by assumption it has finite order.  We conclude using 
Lemmas~\ref{bb} and~\ref{nb}.

For the converse, suppose that every nef divisor on $X$ is semiample.  Let $N$ be a nef divisor on $X$ 
such that $K_X \cdot N = 0$.  By Lemma~\ref{tosto} we deduce that $N_E$ has finite order.  Therefore 
$\Gamma$ is finite, being a finitely generated abelian group all of whose elements have finite order; the 
result follows.
\end{proof}

\begin{example}
Let $X$ be a smooth projective rational surface such that a positive multiple of the anticanonical divisor is effective.

Suppose that the anticanonical divisor of $X$ is not linearly equivalent to an effective divisor, and let $e$ be the least among the positive integers $n$ such that the linear system $|{-nK_X}|$ is not empty.  In view of Theorem~\ref{nof}, every nef divisor on $X$ is semiample.  We first show that if $E$ is an element of the linear system $|{-eK_X}|$, then the cohomology group ${\rm H}^1(E,\mathcal{O}_E)$ vanishes.  Indeed, from the long exact cohomology sequence associated to 
\[
0 \longrightarrow \mathcal{O}_X (eK_X) \longrightarrow \mathcal{O}_X \longrightarrow \mathcal{O}_E \longrightarrow 0
\]
we deduce that the group ${\rm H}^1(E,\mathcal{O}_E)$ is a subgroup of 
\[
{\rm H}^2 \bigl( X,\mathcal{O}_X(eK_X) \bigr) \simeq {\rm H}^0 \bigl( X,\mathcal{O}_X(-(e-1)K_X) \bigr),
\]
which vanishes by the definition of $e$.  
Next, if $N$ is a nef divisor on $X$ such that $N \cdot K_X = 0$, then it follows that the restriction of $N$ to $E$ is trivial, since the divisor $N$ has degree zero on each irreducible component of $E$.  We see that in this case ``there is no space'' for nef divisors that are not semiample: the subgroup of $\Pic(E)$ generated by the classes of the restrictions of the nef divisors to $E$ is trivial.

Suppose that the anticanonical divisor on $X$ is linearly equivalent to an effective divisor, let $E$ be an element of the anticanonical linear system on $X$, and let $\lE \subset \Pic (X)$ be the sublattice generated by the irreducible components of $E$.  We analyze the possibilities for the group $\Gamma$ in terms of the anticanonical Iitaka dimension of $X$ and of the signature of the quadratic form on $\lE$.

The anticanonical Iitaka dimension of $X$ equals two if and only if the quadratic form $\lE$ represents positive values.  In this case the anticanonical divisor of $X$ is big and the group $\Gamma$ is trivial, since a nef divisor that is orthogonal to a big divisor vanishes.  The class of rational surfaces with big anticanonical divisor is analyzed in~\cite{tvv}.

If the anticanonical Iitaka dimension of $X$ equals $1$, then the quadratic form $\lE$ is negative semidefinite and it is not definite.  Let $M$ denote the moving part of a positive multiple of the anticanonical divisor on $X$ with positive dimensional linear system.  A nef divisor on $X$ orthogonal to the anticanonical divisor must be orthogonal to $M$, and, by the Hodge Index Theorem, it is therefore proportional to $M$.  We deduce that $\Gamma$ is cyclic and torsion.  The class of rational surface with anticanonical Iitaka dimension $1$ is analyzed in~\cite{al}.

If the anticanonical Iitaka dimension of $X$ equals $0$, then the quadratic form $\lE$ is negative semidefinite.  If the quadratic form on $\lE$ is not definite, then there is an effective nef divisor $F$ with support contained in the support of $E$, whose class generates the kernel of the intersection form on $\lE$.  In this case, the group $\Gamma$ is again cyclic, contains the group generated by $F$ with finite index, but $\Gamma$ is not torsion, since otherwise a positive multiple of $F$ would be base point free, contradicting the fact that the anticanonical Iitaka dimension vanishes.  If the quadratic form on $\lE$ is negative definite, then there are easy examples in which the group $\Gamma$ is infinite and others in which it is finite.
We analyze one such example after proving the following result.
\end{example}

\begin{prop} \label{neg-def}
Let $X$ and $\lE$ be as above.
Suppose that the quadratic form on $\lE$ is negative definite. Then every nef divisor of $X$ is semiample if and only if, for every $[D] \in \lE^{\perp}$, the line bundle $\cO_E(D)$ represents a torsion element of $\Pic(E)$.
\end{prop}

\begin{proof}
The irreducible components of $E$ span an extremal face of the effective cone and therefore there are nef divisors $N$ on $X$ such that if $F$ is an integral curve on $X$ with $N \cdot F = 0$, then $F$ is a component of $E$.  More precisely, the set of nef divisors with the above property is a non-empty open subset of $\lE^\perp $.  In particular, there is a basis of the space $\lE^\perp$ consisting of nef divisors, and we conclude by Theorem~\ref{fatto}.
\end{proof}

\begin{example}
Let $\overline{E} \subset \bP^2$ be a smooth plane cubic curve and let $p_1,\dots,p_r$ be distinct points lying on $\overline{E}$.  Denote by $X$ the blow-up of $\bP^2$ at the points $p_1,\dots,p_r$, by $E$ the strict transform in $X$ of $\overline{E}$.  We denote by $\ell$ a line in $\mathbb{P}^2$ and by $\ell_E$ the effective divisor on $E$ induced by $\ell$; note that the degree of the divisor $\ell_E$ equals $3$.

If $r\leq 8$, then the quadratic form on $\lE$ is positive definite and every nef divisor is semiample (Lemma~\ref{bbp}). 

In case $r=9$ the intersection form on $\lE$ is negative semidefinite and is not definite since $E^2=0$.  We observe that $E$, which is nef, is semiample if and only if the line bundle $\mathcal{O}_E(E) \simeq \mathcal{O}_E \bigl( 3 \ell_E - (p_1 + \cdots + p_9) \bigr)$ is torsion in $\Pic(E)$ (Corollary~\ref{ii}).

Consider now the case $r\geq 10$, that is when
the quadratic form on $\lE$ is negative definite.
Let $E_1 , \ldots , E_r$ be the exceptional divisors of the blow up morphism $X \to \mathbb{P}^2$.
The classes of $E_1-E_r$, $E_2 - E_r$, \ldots , $E_{r-1}-E_r$ and $\ell - 3E_r$ span $\lE^{\perp}$. 
Moreover these classes restrict to the classes of $p_1-p_r$, $p_2-p_r$, \ldots , $p_{r-1}-p_r$ and $\ell_E-3p_r$ on $\Pic(E)$. Thus we conclude by 
Proposition~\ref{neg-def} that any nef divisor on $X$ 
is semiample if and only if each of the above classes is torsion in $\Pic(E)$.
\end{example}

\section{Examples of rational Mori dream surfaces with $k(-K_X)=0$} \label{seexa}

We begin by recalling the following definition given by Dolgachev and Zhang in~\cite{DoZa}.
\begin{definition}
A {\em Coble surface} $X$ is a smooth projective rational surface with $-K_X$ not effective and $-2K_X$ effective.
\end{definition}
In this section we provide an example of a Coble surface $\tilde{Y}$ with finitely generated Cox ring and $k(-K_{\tilde{Y}})=0$.  Here we consider the following construction.  Inside the moduli space of Enriques surfaces there are exactly two $1$-dimensional families of Enriques surfaces whose general element has finite automorphism group. These have been classified in~\cite{Ko}. To describe the general element $Y_t$ of one such family, we recall a construction given in~\cite{Ko}*{\S3}.  Let $\phi$ be the involution of $\bP^1\times\bP^1$ defined by
\[
  \phi([x_0,x_1],[y_0,y_1])
  =
  ([x_0,-x_1],[y_0,-y_1]).
\]
Let $\bigl\{ C_t \bigr\}_{t \in \bP^1}$ be the pencil of curves of degree $(2,2)$ in $\bP^1\times\bP^1$ given by 
\[
 C_t := V \bigl( (2x_0^2-x_1^2)(y_0^2-y_1^2)
 +
 (2ty_0^2+(1-2t)y_1^2)(x_1^2-x_0^2) \bigr)
\]
and define 
\[
L_1 := V(x_0-x_1) , \quad L_2 := V(x_0+x_1) , \quad L_3 := V(y_0-y_1) , \quad L_4 := V(y_0+y_1).
\]
An elementary calculation shows that $C_t$ is smooth irreducible for $t \neq 1 , \frac{1}{2} , \frac{3}{2} , \infty$. Moreover, the curve $C_t$ has an ordinary double point, for $t \in \{ \frac{1}{2} , \frac{3}{2} \}$, it is the union of two irreducible curves of degree $(1,1)$, for $t=1$, and is the union of $L_1,\dots,L_4$ for $t=\infty$. Note that the base locus of the pencil $C_t$ consists of the four points 
\[
\bigl( [1,1],[1,1] \bigr) , \quad \bigl( [1,1],[1,-1] \bigr) , \quad \bigl( [1,-1],[1,1] \bigr) , \quad \bigl( [1,-1],[1,-1] \bigr), 
\]
and also note that each of these points is contained in two of the curves $L_1,\ldots,L_4$.  For each $s \in \bP^1$ the birational map $\pi_s \colon S_s \to \bP^1 \times \bP^1$ is obtained by first blowing up the four base points of the pencil $C_t$ and then blowing up the $12$ points of intersections of any exceptional divisor with the strict transform of $C_s$ and of $L_1 , \ldots , L_4$.  Thus the surface $S_s$ is rational with Picard group of rank $18$.

For each $t \in \bP^1$, the reducible curve $B_t \subset \bP^1 \times \bP^1$ defined by 
\[
 B_t := C_t+L_1+L_2+L_3+L_4
\]
is a divisor of degree $(4,4)$, invariant with respect to the involution $\phi$.  Dropping the subscript $t$ from the morphisms to simplify the notation, this leads us to the diagram
\[
 \xymatrix{
 X_t\ar[r]^-{\varphi}_{2:1}\ar[d]_-{\psi}^-{2:1} 
 & S_t\ar[r]^-{\pi} & \bP^1\times\bP^1\\
 Y_t
 }
\]
where the morphism $\pi$ is birational, the surface $S_t$ is rational with Picard group of rank $18$, the morphism $\varphi$ is a double cover, the branch locus of $\varphi$ consists of the strict transform $B'_t$ of $B_t$ together with a union $\Gamma$ of disjoint curves in the exceptional locus of $\pi$, and finally $\psi$ is the double cover $X_t \to Y_t = X_t / \langle \sigma\rangle$, where $\sigma \in \Aut(X_t)$ is the involution induced by $\phi$.

Kondo proves that $Y_t$ is an Enriques surface for all $t$ different from $1 , \frac{1}{2} , \frac{3}{2} , \infty$. Moreover he proves that $\Aut(Y_t)$ is a finite group for any value of $t$ (see~\cite{Ko}*{\S3 and Corollary~5.7} for precise statements).  Let $t_0 \in \bP^1$ be such that $C'_{t_0}$ is irreducible and has an ordinary double point and let $X_0 := X_{t_0}$ be the corresponding surface. Then $X_0$ has a singularity of type $\mathbf{A}_1$ at a point $p$. Since $p$ is the only singular point of $X_0$, it must be stable with respect to the involution $\sigma$, and thus $Y_0$ is singular at one point $q$ as well. This leads to the commutative diagram
\[
 \xymatrix{
 \tilde{X}\ar[r]^-{\pi_p}
 \ar[d]_-{\tilde{\psi}} & X_0\ar[d]^-{\psi}\\
 \tilde{Y}\ar[r]_-{\pi_q} & Y_0
 }
\]
where $\pi_p$ and $\pi_q$ are minimal resolution of singularities, $\tilde{X}$ is a K3-surface and $\tilde{\psi}$ is a double cover branched along the exceptional divisor $E$ of $\pi_q$.

\begin{prop}
The surface $\tilde{Y}$ is a Coble surface.
\end{prop}

\begin{proof}
Let $\sigma\in\Aut(X_0)$ be the involution induced by $\phi$.  As we noted before, $\sigma(p)=p$ so that $\sigma$ lifts to an involution $\tilde{\sigma}$ of $\tilde{X}$. The involution $\tilde{\sigma}$ is non-symplectic since $\sigma$ is non-symplectic and that the exceptional divisor $R:=\pi_p^{-1}(p)$ is fixed by this involution. Thus the quotient surface $\tilde{Y}=\tilde{X}/\langle\tilde{\sigma}\rangle$ is a smooth projective rational surface.  Observe that the branch divisor of $\tilde{\psi}$ is contained in $|{-2K_{\tilde{Y}}}|$, and in particular $-2K_{\tilde{Y}}$ is effective. Moreover, since $\psi$ does not ramify along a divisor, the branch divisor of $\tilde{\psi}$ is exactly $E$. Since $X_0$ has an $\mathbf{A}_1$-singularity at $p$, it follows that $F=\pi_p^{-1}(p)$ is a $(-2)$-curve on $\tilde{X}$. In particular, the curve $E=\tilde{\psi}(F)$ is irreducible and reduced with $E^2<0$. Thus $-K_{\tilde{Y}}$ cannot be effective, since otherwise $E \in |{-2K_{\tilde{Y}}}|$ would be non-reduced.
\end{proof}
We are now ready to prove the following.

\begin{theorem}
The Cox ring of $\tilde{Y}$ is finitely generated.
\end{theorem}

\begin{proof}
We want to use~\cite{Ko}*{Corollary~5.7} to prove that the group $\Aut(Y_0)$ is finite, even though $Y_0$ is a rational surface and not an Enriques surface.  The reason we can still apply Kondo's result is that the N\'eron-Severi lattice of $Y_0$ is isomorphic to the N\'eron-Severi lattice of an Enriques surface, and moreover the minimal resolution of the limit $X_0$ of the K3-surfaces $X_t$ is still a K3-surface and therefore it still corresponds to a point in the period space.  Moreover $Y_0$ is Kawamata log terminal, klt for short, since 
\[
K_{\tilde{Y}} = \pi_q^*K_{Y_0}-\frac{1}{2}E
\]
and $\pi_q$ is a resolution of singularities of $Y_0$. Observe that indeed $Y_0$ is a klt Calabi-Yau surface since its canonical divisor is numerically trivial. Let $\Delta:=\frac{1}{2}E$. Since the pair $(\tilde{Y},\Delta)$ is the terminal model of $Y_0$, it follows that $\Aut(Y_0)=\Aut(\tilde{Y},\Delta)$. In particular the last group is finite so that the image of the map $\Aut(X,\Delta)\to {\rm GL} \bigl( \Pic_\bQ(X) \bigr)$ is finite as well. Hence the Cox ring of $\tilde{Y}$ is finitely generated by~\cite{To}*{Corollary~5.1}.
\end{proof}

\begin{remark}
In general a Coble surface is not a Mori dream space. To see this let $p_1,\dots,p_9$ denote the nine intersection points of two general plane cubic curves $C_1$ and $C_2$.  Let $q\in C_1$ be such that the divisor class of $p_9-q$ is a $2$-torsion point of $\Pic(C_1)$. Then the blow up $Z$ of the plane at $p_1,\dots,p_8,q$ admits an elliptic fibration defined by $|{-2K_Z}|$. Due to the generality assumption on $C_1$ and $C_2$, it is easy to see that the linear system $|{-2K_Z}|$ does not contain reducible elements. Equivalently $Z$ does not contain $(-2)$-curves, so that it is not a Mori dream space by~\cite{al}.
By an Euler characteristic calculation the fibration induced by $|{-2K_Z}|$ contains $12$ nodal curves. The blow up $Y$ of $Z$ at one of these nodes $p$ is a Coble surface since $|{-K_Y}|$ is empty, but $|{-2K_Y}|$ contains the strict transform of the fiber through $p$. Moreover $Y$ is not a Mori dream surface since $Z$ is not a Mori dream surface.
\end{remark}

In the following remark, we obtain a description of the Coble surface $\tilde{Y}$ as an iterated blow up of $\mathbb{P}^2$.

\begin{remark}
The Coble surface $\tilde{Y}$ constructed above as the minimal resolution of a limit of a family of Mori dream Enriques surfaces contains (at least) the configuration of $(-2)$-curves of the general element of such family.  The intersection graph of this surface is given in Figure~\ref{tily} (see~\cite{Ko}*{Example~1}), where double edges mean that the corresponding curves have intersection equal to $2$.

\begin{figure}[h]
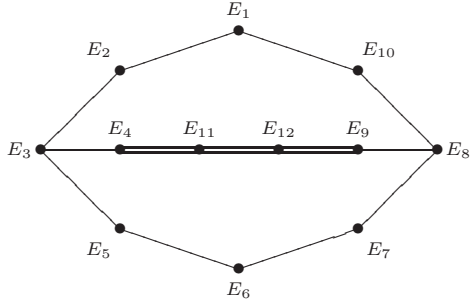

\begin{minipage}{200pt}
\[
\xygraph {[] !~:{@{=}}
!{<0pt,0pt>;<15pt,0pt>:}
!{\save +<0pt,8pt>*\txt{$\scriptstyle E_1$}  \restore}
{\bullet} [drrr]
!{\save +<8pt,8pt>*\txt{$\scriptstyle E_{10}$}  \restore}
{\bullet} [ddrr]
!{\save +<8pt,0pt>*\txt{$\scriptstyle E_8$}  \restore}
{\bullet} [ddll] 
!{\save +<8pt,-8pt>*\txt{$\scriptstyle E_7$}  \restore}
{\bullet} [dlll]
!{\save +<0pt,-8pt>*\txt{$\scriptstyle E_6$}  \restore}
{\bullet} [ulll]
!{\save +<-8pt,-8pt>*\txt{$\scriptstyle E_5$}  \restore}
{\bullet} [uull]
!{\save +<-8pt,0pt>*\txt{$\scriptstyle E_3$}  \restore}
{\bullet} [uurr] 
!{\save +<-8pt,8pt>*\txt{$\scriptstyle E_2$}  \restore}
{\bullet} [dd] 
!{\save +<0pt,8pt>*\txt{$\scriptstyle E_4$}  \restore}
{\bullet} [rr] 
!{\save +<0pt,8pt>*\txt{$\scriptstyle E_{11}$}  \restore}
{\bullet} [rr] 
!{\save +<0pt,8pt>*\txt{$\scriptstyle E_{12}$}  \restore}
{\bullet} [rr] 
!{\save +<0pt,8pt>*\txt{$\scriptstyle E_9$}  \restore}
{\bullet} [rl]
- [rr] - [uull] - [ulll] - [dlll] - [ddll] - [ddrr] -[drrr] -[urrr] - [uurr] [ll] 
: [llllll]
- [ll]
}
\]
\end{minipage} 
\caption{Configuration of $(-2)$-curves of $\tilde{Y}$.} \label{tily}
\end{figure}

It is not difficult to see that $F:=E_1+E_2+E_3+E_5+E_6+E_7+E_8+E_{10}$ has $F^2=0$ and $|2F|$ is an elliptic pencil on $\tilde{Y}$ with two reducible fibers one of type $I_8$ and one of type $I_2$. Thus $\tilde{Y}$ is the blow-up of a Mori dream surface $Z$ with Picard group of rank $10$, admitting an elliptic fibration whose dual graph of singular fibers contains $I_8$ and $I_2$. The surface $Z'$ defined by the jacobian fibration of that induced by $|{-2K_Z}|$ has finite Mordell-Weil group by~\cite{al}. Moreover $Z'$ is unique up to isomorphism and its unique elliptic fibration, given by $|{-K_{Z'}}|$, admits exactly four singular fibers of type $I_8$, $I_2$, $I_1$, $I_1$, by~\cite{Du}*{Proposition 9.2.19}. This is the same configuration of singular fibers of the elliptic pencil defined by the linear system $|{-2K_Z}|$.  Explicitly, let $x,y,z$ be homogeneous coordinates on $\mathbb{P}^2$ and let 
\begin{eqnarray*}
& c_0 := 4 x^2 (x^2 - yz)^2 \\[5pt]
&{\textrm{and}} \\[5pt]
& c_\infty := (2x^2y - x^2z + y^3 - 2y^2z + yz^2) (2x^2y + x^2z - y^3 - 2y^2z - yz^2).
\end{eqnarray*}
The surface $Z$ is the rational elliptic surface associated to the pencil generated by the two plane sextics $C_0 = V(c_0)$ and $C_\infty = V(c_\infty)$.  The fiber corresponding to the curve $C_0$ has multiplicity two and is of type $I_8$.  The fiber corresponding to the curve $C_\infty$ is reduced and is of type $I_2$.  The only two remaining singular fibers are both of type $I_1$ and correspond to the plane sextics $V(c_0 + c_\infty)$ and $V(c_0 - c_\infty)$.  To construct the jacobian surface $Z'$ associated to $Z$, let 
\[
d_0 := 4x(yz-x^2) \quad \quad {\textrm{and}} \quad \quad d_\infty := z(y + z)^2 - x^2(y + 2z);
\]
the surface $Z'$ is the rational elliptic surface determined by the pencil generated by the plane cubics $D_0 = V(d_0)$ and $D_\infty = V(d_\infty)$.  As before, the fiber corresponding to $D_0$ is of type $I_8$; the fiber corresponding to $D_\infty$ is of type $I_2$; the only remaining singular fibers are of type $I_1$ and correspond to the cubics $V(d_0+d_\infty)$ and $V(d_0-d_\infty)$.  Observe that the singularity of the pencil defining $Z'$ at the point $[0,1,-1]$ is resolved only after two successive blow ups.  To obtain $Z$ from $Z'$ it suffices to contract the exceptional curve lying above the point $[0,1,-1]$ (the one introduced by the second blow up) and to blow up the point $[0,1,1]$.  Finally, the Coble surface $\tilde{Y}$ is obtained from $Z$ by blowing up the singular point of one of the two $I_1$ fibers.  The two surfaces obtained by the choice of the last blown up point are isomorphic: the substitution $([x,y,z],t) \mapsto ([ix,y,-z],-t)$ in the pencil $c_0+tc_\infty$ defining $Z$ determines an automorphism of order 4 on $Z$ exchanging the two $I_1$ fibers.  
\end{remark}

\end{document}